\documentclass[12pt]{report}
\usepackage{amsmath}
\usepackage{amsfonts}
\usepackage{amssymb}
\usepackage{mathrsfs}
\usepackage{graphicx}
\usepackage{cancel}
\usepackage{hyperref}
\usepackage[english]{babel}
\addto{\captionsenglish}{%
  
}
\DeclareMathAlphabet{\mathpzc}{OT1}{pzc}{m}{it}

\newtheorem{theorem}{THEOREM}
\newtheorem{definition}{DEFINITION}
\newtheorem{lemma}{LEMMA}

\begin{document}
\pagenumbering{roman}
\begin{flushleft}
JOURNAL OF COMBINATORIAL THEORY \textbf{(B) 22,} 91-96(1977)
\end{flushleft}
$$$$$$$$
\begin{center}
\textbf{\large Circuit Preserving Edge Maps $$$$
JON HENRY SANDERS
}
\end{center}
\begin{center}
\textit{jon\_sanders@partech.com}
\end{center}
\begin{center}
AND
\end{center}
\begin{center}
\textbf{\large DAVID SANDERS
}\end{center}

\begin{center}
\textit{davidhsanders@earthlink.net}
\end{center}
$$$$
\begin{center}
\textit{communicated by A.J. Hoffman}
\end{center}
$$$$$$$$
\begin{center}
Received March 12, 1975
\end{center}$$$$

It is proved than any one-to-one edge map $f$ from a 3-connected graph $G$ onto a graph $G^\prime$,$G$ anf $G^\prime$ possibly infinite, satisfying $f(C)$ is a circuit in $G^\prime$ whenever $C$ is a circuit in $G$ is induced by a vertex isomorphism. This generalizes a result of Whitney which hypothesizes $f(C)$ is a circuit in $G^\prime$ if and only if $C$ is a circuit in $G$.
\section*{\small 
\begin{center}
1. INTRODUCTION
\end{center}
}
\pagenumbering{arabic}
\setcounter{page}{1}
\noindent
\par In 1932, Whitney proved [3] that every circuit isomorphism (one-to-one onto edge map $f$ such that $C$ is a circuit if and only if $f(C$) is a circuit) between two 3-connected graphs is induced by a vertex isomorphism. The following year Whitney observed [4] that this result could be strengthened by hypothesizing the 3-connectivity of only one of the graphs. It is necessary to also assume the other graph has no isolated vertices. In 1966, Jung pointed out[1] that Whitney's result also holds for infinite graphs.
\par In this paper, we further generalize Whitney's result by proving that any circuit injection $f$ (a one-to-one edge map such that if $C$ is a circuit then $f(C)$ is a circuit) from a 3-connected graph $G$ onto a graph $G^\prime$ is induced by a vertex isomorphism. Throughout we will understand the terminology that $f$ is a circuit injection from $G$ onto $G^\prime$ to preclude the possibility of $G^\prime$ having isolated vertices. The term graph refers to undirected graphs, finite or infinite , without loops or multiple edges.
\par We note that a circuit injection $f:G\rightarrow G^\prime$ where $G$ is 2-connected is not necessarily a vertex or circuit isomorphism no matter what connectivity $n$ is assumed for $G^\prime$ as illustrated by the following example. For any prime $p>2$ let $G$ be the graph consisting of $p$ paths $P_i, i\in Z_p$(where $Z_p$ is the integers modulo $p$), each path having the same two endpoints but otherwise mutually disjoint, and each $P_i$ consisting of $p$ edges $e_{i\cdot j}$,$j\in Z_p$. Let $G^\prime$ be the complete bipartite graph on the vertex sets$\left\{b_i:i\in Z_p\right\}$ and $\{c_i:i\in Z_p\};$ and define the edge map $f:G\rightarrow G^\prime$ by $f(e_{i\cdot j})=(b_j,c_{i+j})$ where $i\in Z_p, j\in Z_p$ . Then $G$ is 2-connected, $G^\prime$ is $p$-connected and it can be checked that $f(C)$ is a circuit whenever $C$ is a circuit.
\newpage
\section*{\small{
\begin{center}
2. THEOREMS AND PROOFS
\end{center}
}}
\noindent
\par Our principal result is Theorem 6 whose  Proof consists of the application of Theorems 1 through 5.

\begin{theorem}
Let $G$ and $G^\prime$ be graphs without isolated vertices, $G$ without isolated edges, and $g:G\rightarrow G^\prime$ is a one-to-one map of the edges of $G$ onto the edges of $G^\prime$ such that for each vertex $v$ of $G$ the star subgraph $S(v)$ is mapped by $g$ onto the star subgrapgh $S(v^\prime)$ for some vertex $v^\prime$ of $G^\prime$. Then $g$ is induced by a vertex isomorphism $\lambda$.
\end{theorem}
\textit{Proof.} ~For each vertex $v$ of $G$ let $\lambda(v)=v^\prime$ be a vertex such that $g(S(v))=S(v^\prime)$.
It can be verified that $v^\prime$ is then uniquely determined, but this is not necessary. To see that $\lambda$ is one-to-one note that if $\lambda(u)=\lambda(v)$  then $S(\lambda(u))=S(\lambda(v)),$ thus $g(S(u))=g(S(v)),$ which implies 
$S(u)=S(v),$ which implies $u=v$, edge$(u,v)$ is isolated, or $u$ and $v$ are isolated vertices. To see that $\lambda$ is onto, given any vertex $w$ of $G^\prime$ let $e$ be an edge incident to $w$ and then using the definition of $\lambda$ and that $\lambda$ is one-to-one it is seen that $\lambda$ must map one of the vertices of $g^{-1}(e)$ into $w$. To see that $\lambda$ induces $g$, we observe that there exists an edge $(\lambda(u),\lambda(v))$ in $G^\prime$ if and only if 
$S(\lambda(u))\cap S(\lambda(v))\neq\phi $ if and only if $g^{-1}(S(\lambda(u))\cap S(\lambda(v)))\neq\phi $  if and only if 
$g^{-1}(S(\lambda(u)))\cap g^{-1}(S(\lambda(v)))\neq \phi$ if and only if $S(u)\cap S(v)\neq\phi $ if and only if there exists an edge $(u,v)$ in $G$.

\begin{lemma}
Let $a,b,c$ be three distinct vertices of a 2-connected graph $G$. Then there exists a circuit $C$ containing $a$ and $b$ and a path $P(c,t)$ where $t$ is a vertex on $C$ different from $a$ and $b$ and no other vertex of $P(c,t)$ is on $C$. We allow the possibility $c=t$ and $P(c,t)=\phi$.
\end{lemma}
\textit{Proof.}~~
Take any circuit containing $a$ and $b$. if $c$ is on $C$ then we have the case with $P(c,t)=\phi$. If $c$ is not on $C$ choose any vertex $v$ of $C, v\neq a, v\neq b$ and let $C_1=P_1(c,v)\cup P_2(c,v)$ 
be a circuit through $c$ and $v$. Let $t_1$ and $t_2$ be the first vertices of $P_1(c,v)$ respectively $P_2(c,v)$ which lie on $C$. If 
$\left\{t_1,t_2\right\}=\left\{a,b\right\}$ 
then $C_1$ is a circuit containing $a,b,c$ and again we have the case with $P(c,t)=\phi$. Otherwise at least one of the $t_i$ is different from $a$ and $b$ and the corresponding $P_i(c,t_i)$ with $C$ are desired path and  circuit.

\begin{lemma}
Let $f$ be a circuit injection from $G$ onto $G^\prime$, $G$ 3-connected, and $S(v)$ a star subgraph of $G$. Then $f(S(v))$ is either a star subgraph of $G^\prime$ or an independent (i.e., pairwise nonadjacent) set of edges.
\end{lemma}

\textit{Proof.}~~
If $f(S(v))$ is not an independent set of edges then there are two edges $e_1=(a_1,v)$ and $a_2,v$ of $S(v)$ with $f(e_1)$ and $f(e_2)$ adjacent in $G^\prime$ at some vertex $w$. Suppose some other edge $e_3=(a_3,v)$ of $S(v)$ does not have its image $f(e_3)$ incident to $w$. Since $G-v$ is 2-connected, by Lemma 1 there is a circuit $C=P_1(a_1,a_3)\cup P_2(a_1,a_3)$
and a path $P(a_2,t)$ with no vertex on $C$ except t. $ C_1=P_1\cup\{e_1,e_3\}$ is a circuit in $G$ so $f(C_1)=f(P_1)\cup\{f(e_1),f(e_3)\}$ is a circuit in $G^\prime$. By hypothesis $f(C_1)$ passes through $w$ and $f(e_2)$ does not. So some edge $f(p_1)$ of $f(P_1)$ must be incident to $w$. Similarly some edge $f(P_1)$ of $f(p_2)$ must be incident to  $w$. We derive a contadiction to $f(p_1),f(p_2),f(e_2)$ each incident to $w$ by finding a circuit in $G$ containing $p_1,p_2$, and $e_2$. Since $t$ lies on $C$, we have $t$ on $P_1$ or $P_2$. Suppose without loss of generality $t$ lies on $P_1$ so that we may write $P_1(a_1,a_3)=P_1(a_1,t)\cup P_1(t,a_3)$. If $p_1$ is on $P_2(a_1,t)$ then the circuit $P_1(a_1,t)\cup P(a_2,t)\cup\{e_2,e_3\}\cup P_2(a_1,a_3)$ contains $p_1,p_2$ and $e_2$. If $p_1$ is  on $P_1(t,a_3)$, then the desired circuit is $P_1(t,a_3)\cup P(a_2,t)\cup\{e_1,e_2\}\cup P_2(a_1,a_3)$.
\par Thus we have shown that if $f(S(v))$ is not an independent set of edges $f(S(v))$ is a subset of a star subgraph $S(w)$ of $G^\prime$. To finish the proof suppose there were some edge $f(e_4)$ at $w$ with $e_4\notin S(v).$ Pick any edge $e$ of $S(v)$ and a circuit $C^\prime$ in $G$ contaning $e$ and $e_4$. $C^\prime$ must contain another edge $e^\prime$ of $S(v)$ but then we have the contradiction that $f(C^\prime)$ cannot be a circuit because $f(e)\cdot f(e^\prime)$, and $f(e_4)$ are each incident at $w$.

\begin{theorem}
Let $f$ be a circuit injection from $G^\prime$ onto $G$ 3-connected, and $S(w)$ a star subgraph of $G_\prime$. Then $f^{-1}(S(w))$ is either a star subgraph of $G$ or an independent set of edges.

\end{theorem}

 \textit{Proof.}~~~~
If $f^{-1}(S(w))$ is not an independent set of edges, then there exist $e_1$ and $e_2\in f^{-1}(S(w))$ such that $e_1$ and $e_2$ have common vertex $v$. By lemma 2, $f(S(v))$ is either an independent set or a star subgraph of $G^\prime$.The former case is ruled out since $f(e_1)$ and $f(e_2)$ are adjacent at $w$. Thus $f(S(v))=S(w^\prime)$ for some vertex $w^\prime$ of $G^\prime$. But since $\{f(e_1),f(e_2)\}\subset S(w)\cap S(w^\prime)$ we have $w=w^\prime.$
Thus $f(S(v))=S(w),$ hence $f^{-1}(S(w))=S(v)$.

\begin{theorem}
Let $f$ be a circuit injection from  $G$ onto $G^\prime, G$ 2-connected, and $S=S(v)$ a star subgraph of $G^\prime.$ Then $G=G_1\cup G_2\cup f^{-1}(S)$, where $G_1$ and $G_2$ are connected components of $G-f^{-1}(S)$. (with $G-f^{-1}(S)$) denoting the subgraph of $G$ containing the same vertices as $G$ but only those edges of $G$ not in $f^{-1}(S)$ and where each edge of $f^{-1}(S)$ has one vertex in $G_1$ and one vertex in $G_2.$
\end{theorem}
 
\textit{Proof.}~~~~
Let $G_\alpha,\alpha\in I$ be the connected components of $G-f^{-1}(S)$. Each edge $e=(a,b)\in f^{-1}(S)$ cannot have both vertices $a,b$ in the same connected component $G_\alpha$, for otherwise there would exist a path $P(a,b)\subset G_\alpha$ , a circuit $C=\left\{e\right\}\cup P(a,b)$ and therefore a circuit $f(C)$ containing only one edge $f(e)$ of $S(v),$ an impossibility. It remains only to show $|I|=2.$ From the preceding, $|I|>1$, so assume $|I|\geq 3$. Take any three connected components $G_1,G_2,G_3$ of $G-f^{-1}(S).$ If   there were edges $e_{12}=(a_1,a_2).e_{23}=(b_2,b_3), e_{31}=(c_3,c_1)$ of $f^{-1}(S)$ joining  $G_1$  to $G_2$, $G_2$ to $G_3$, $G_3$ to $G_1$ , respectively, there would be a circuit $C_1$ in $G$ consisting of $\{e_{12},e_{23},e_{31}\}$ and paths $P(c_1,a_1)$ in $G_1, P(a_2,b_2)$ in $G_2$, and $P(b_3,c_3)$ in $G_3$. Then we have the contradiction that there is a circuit $f(C_1)$ in $G^\prime$ containing three edges $f(e_{12}),f(e_{23}),$ and $f(e_{31})$ of $S(v)$. So at least two of the components, say $G_1$ and $G_2$, are not joined by any edge of $f^{-1}(S).$ Choose a vertex $v_1$ in $G_1$ and a vertex $v_2$ in $G_2$. Since $G$ is 2-connected there is a circuit $C_2$ in $G$ containing $v_1$ and $v_2,C_2$ consisting of two paths $P_1(v_1,v_2)$ and $P_2(v_1,v_2)$ having only $v_1$ and $v_2$ in common. Because no edge of $f^{-1}(S)$ joins $G_1$ and $G_2, P_1$ and $P_2$ each contain two edges of $f^{-1}(S)$. But then we have the contradiction that $f(C_2)$ contains four or more edges of $S(v)$. This $|I|=2$ and the Proof is complete.
\begin{definition}
Let $G$ be a graph consisting of two vertex disjoint circuits $A$ and $B$, two  edges $e_1=(a_1,b_1), e_2=(a_2,b_2)$ and a path $P(a_3,b_3)$ vertex disjoint except for $a_3$ and $b_3$ from $A$ and $B$, where $a_1,a_2,a_3$ are distinct vertices of $A$ and $b_1,b_2,b_3$ are distinct vertices of $B$. Let $e_3$ be an arbitrary edge of $P(a_3,b_3)$. We say $G$ is a graph of type $X$ with connectors $e_1,e_2,$ and $e_3$.
\end{definition}
\begin{theorem}
Let $G$ be 3-connected and let $A=\{e_1,e_2,\cdots,e_n\}$ be a set of independent edges of $G$ such that $G$ - $A$ has two connected components $G_1$ and $G_2$ and each edge of $A$ has one vertex in $G_1$ and one vertex in $G_2$. Then either $G$ has a subgraph of type $X$ with three connectors from $A$ or there exists a circuit containing at least four distinct edge in $A$.
\end{theorem}
\textit{ Proof.}~~ We consider two cases.\\
\newline
\textit{Case 1.}~~ $G_1$ and $G_2$ are both 2-connected. By the 3-connectivity of $G$ there must be at least three edges in $A$, $ e_1=(a_1,b_1), e_2=(a_2,b_2),$ and $e_3=(a_3,b_3)$ with the $a's$ distinct and in $G_1$, the $b's$ distinct and in $G_2$. By Lemma 1 there exist a circuit $C_1$ containing $a_1$ and $a_2$ and a path $P_1(a_3,t)$ having  no vertex in common with $C_1$ except $t$ which is different from $a_1,a_2.$ Similarly, there is a circuit $C_2$ containing $b_1$ and $b_2$ and a path $P_2(b_3,t^\prime)$ vertex disjoint from $C_2$ except for $t^\prime\neq b_1, b_2.$ Then $C_1,C_2,\{e_1,e_2\},$ and $P_1(a_3,t)\cup \{e_3\}\cup P_2(b_3,t^\prime)$ constitute a subgraph of type $X$ with connectors $e_1, e_2$ and $e_3$.\\
\newline
\textit{Case 2.}~~ $G_1$ and $G_2$ are not both 2-connected. Then at least one of $G_1$ and $G_2,$ say $G_1$ has a cutpoint $v$. Choose vertices $a$ and $b$ in different components of $G_1-v$. By the 3-connectivity of $G$ there are two paths $P_1(a,b)$ and $P_2(a,b)$ in $G$- $v$ having only $a$ and $b$ in common, and each of these paths must have at least two edges of $A$. This gives a circuit containing at least four distinct edges of $A$.
\begin{theorem}
If $G$ is a graph of type $X$ with connectors $e_1,e_2,e_3\in P(a_3,b_3)$ and $f$ is a circuit injection from $G$ onto $G^\prime$, then $f(e_1)$ and $f(e_2)$ do not have a common vertex.
\end{theorem}

\textit{Proof.} 
For any edge, path, or circuit $P$ of $G$ let $P^\prime=f(P).$ Suppose $f(e_1)$ and $f(e_2)$ have a common vertex so we may write $e_1^\prime=(v_1,v_0)$ and $e^\prime_2=(v_2,v_0)$. In the notation of Definition 1 we may also write $A=P(a_1,a_2)\cup P(a_2,a_3)\cup P(a_3,a_1)$ and $B=P(b_1,b_2)\cup P(b_2,b_3)\cup P(b_3,b_1).$ Since $\{e_1,e_2\}\cup P(a_1,a_2)\cup P(b_1,b_2)$ is a circuit in $G$, $\{e_1^\prime,e_2^\prime\}\cup P^\prime(a_1,a_2)\cup P^\prime(b_1,b_2)$ is a circuit in $G^\prime$. Thus the edges of $P^\prime(a_1,a_2)\cup P^\prime(b_1,b_2)$ form a path $P(v_1,v_2)$. Let $v\neq v_1, v_2$ be a vertex in $G^\prime$ where an edge $e_0^\prime$ of $P^\prime(a_1,a_2)$ and an edge of $P^\prime(b_1,b_2)$ meet. $A^\prime$ is a circuit containing $P^\prime(a_1,a_2)$ and disjoint from $P^\prime(b_1,b_2)$. Let $e^\prime$ be an edge of $A^\prime$ at $v$, $e^\prime\neq e_0^\prime$. We have $e^\prime\notin P^\prime(a_1,a_2)$ since otherwise  there would be two edges of $P^\prime(a_1,a_2)$ and an edge of $P^\prime(b_1,b_2)$ incident at $v$ contradicting $P^\prime(a_1,a_2)\cup P^\prime(b_1,b_2)$ is a path. Also, $e^\prime\notin P^\prime(a_2,a_3)$ since otherwise $v$ is a vertex of degree at least 3 in the subgraph $P^\prime(a_1,a_2)\cup P^\prime(a_2,a_3)\cup P^\prime(b_1,b_2)$ which is contained in the circuit  $P^\prime(a_1,a_2)\cup P^\prime(a_2,a_3)\cup P^\prime(a_3,b_3)\cup  P^\prime(b_2,b_3)\cup P^\prime(b_1,b_2)\cup\{e_1\}$. Similarly, $e^1\notin P^\prime(a_3,a_1)$ since otherwise $v$ has degree at least 3 in the subgraph $P^\prime(a_1,a_2)\cup P^\prime(b_1,b_2)\cup P^\prime(a_3,a_1)$ which is contained in the circuit $P^\prime(a_1,a_2)\cup P^\prime(a_3,a_1)\cup P^\prime(a_3,b_3)\cup  P^\prime(b_3,b_1)\cup P^\prime(b_1,b_2)\cup\{e_1\}$. Thus we have a condition to $e^\prime\in A^\prime=P^\prime(a_1,a_2)\cup P^\prime(a_2,a_3)\cup P^\prime(a_3,a_1)$ and the Proof is complete.

\begin{theorem}
If $f$  is a circuit injection from $G$ onto $G^\prime$ where $G$ is 3-connected, then $f$ is induced by a vertex isomorphism.
\end{theorem}
\textit{Proof.}~~~~
We prove  $f$ is induced by a vertex isomorphism by applying Theorem 1 to $f^{-1}$ to show it is induced by a vertex isomorphism. Note Theorem 1 can apply to $f^{-1}$ since $G^\prime$ has no isolated veertices by the assumption that $f$ is onto, and no isolated edges by the fact that any two edges $e_1$ and $e_2$ of $G^\prime$ must lie on some circuit $f(C)$ where $C$ is a circuit containing $f^{-1}(e_1)$ and $f^{-1}(e_2).$ To complete the Proof we must show for any star subgraph $S(v)$ of $G^\prime$ that $f^{-1}(S(v))$ is also a star subgraph. Theorems 2 and 3 tell us the only other possibility for $f^{-1}(S(v))$ is that it is a set of independent edges of $G$ such that  $G-f^{-1}(S(v))$ consists of two connected components $G_1$ and $G_2$ with each edge of $f^{-1}(S(v))$ having one vertex in $G_1$ and one vertex in $G_2.$ But in this event Theorem 4 asserts that either three edges of $f^{-1}(S(v))$ are connectors in a subgraph of $G$ of type $X$ or atleast four edges of $f^{-1}(S(v))$ lie on some circuit $C^\prime$ in $G$. The first situation is ruled out by Theorem 5.\\
The second case is also impossible since it implies $|f(C^\prime)\cap S(v)|\geq 4$ and the theorem is proved.

\newpage
\section*{\small 
\begin{center}
3. GENERALIZATIONS
\end{center}
}
\noindent
\par Possible generalization of Theorem 6 could be attempted by dropping the hypothesis that $f$ is one-to-one. An interesting result of dropping this hypothesis is that the theorem remains true for finite 3-connected graphs, but not for infinite graphs of arbitrarily large connectivity.\\
\par Further generalization could follow the route of assuming $G^\prime$ is not necessarily a graph but a (binary) matroid. Using Tutte's definition of 3-connected for matroids [2], $G$ could also be assumed to be a matriod. The existence of these generalizations will be explored in a following paper.
\\
\begin{center}
\section*{\small ACKNOWLEDGMENT}

\end{center}

\par We express our thanks to the referee for pointing out a simplification in the  Proof of Theorem 4 which eliminates the need to treat infinite graphs separately.

\section*{
\begin{center}
\small REFERENCES
\end{center}
}
\begin{enumerate}
\item H. A. JUNG, Zu einem Isomorphiesatz von H. Whitney f$\ddot{u}$r Graphen, \textit{Math. Ann.} \textbf{164} (1966), 270-271
\item W.T. TUTTE, Connectivity in matroids, \textit{Canad. J. Math.} \textbf{18} (1966), 1301-1324
\item H. WHITNEY, Congruent graphs and the connectivity of graphs, \textit{Amer. J. Math.} \textbf{54} (1932), 150-168
\item H. WHITNEY, 2-isomorphic graphs, \textit{Amer. J. Math.} \textbf{55} (1933), 245-254. 
\end{enumerate}
\end{document}